\documentclass[16pt]{article}
\usepackage{graphicx}
\usepackage{amssymb,verbatim}

\newtheorem{theorem}{Theorem}
\newtheorem{proposition}{Proposition}
\newtheorem{corollary}{Corollary}
\newtheorem{lemma}{Lemma}

\newtheorem{definition}{Definition}

\begin{document}

\title{A proof of  Tait's Conjecture on alternating $-$achiral knots}

\author{Nicola ERMOTTI, Cam Van QUACH HONGLER and Claude WEBER}

\maketitle

\begin{abstract}
In this paper we are interested in symmetries of alternating knots, more precisely in those related to achirality. We call the following statement {\bf Tait's Conjecture on alternating $-$achiral knots}: 

Let $K$ be an alternating $-$achiral knot. Then there exists a minimal projection $\Pi$ of $K$ in $S^2 \subset S^3$ and an involution $\varphi: S^3 \rightarrow S^3$ such that:
\\ 1) $\varphi$ reverses the orientation of $S^3$;
\\ 2) $\varphi(S^2) = S^2$;
\\ 3) $\varphi (\Pi) = \Pi$; 
\\ 4) $\varphi$ has two fixed points on $\Pi$ and hence reverses the orientation of $K$.

The purpose of this paper is to prove this statement. 

For the historical background of the conjecture in Peter Tait's and Mary Haseman's papers see \cite{QW1}. 

\end{abstract}

\section{Introduction}

In this paper we are interested in symmetries of alternating knots, more precisely in those related to achirality. 

\begin{definition}
A knot $K \subset S^3$ is said to be {\bf achiral} if there exists a diffeomorphism $\Phi : S^3 \rightarrow S^3$ such that:
\\
1) $\Phi$ reverses the orientation of $S^3$;
\\
2) $\Phi (K) = K$.
\\
If there exists a diffeomorphism $\Phi$ such that its restriction to $K$ preserves (resp. reverses) the orientation of $K$, the knot is said to be {\bf +achiral} (resp. {\bf -achiral}).
\end{definition}

In his pioneering attempt to classify knots, Peter Tait (see \cite{Tait1} \cite{Tait2} and \cite{Tait3}) considered mostly alternating knot projections. One can see throughout his papers that he was guided by the conviction that most questions about alternating knots can be answered by an effective procedure based on minimal alternating knot projections. Several  unproved principles used by Tait were later called ``Tait's Conjectures" and proved in the eighties and nineties of last century, after Vaughn Jones'  discovery  of his polynomial. On the other hand, it is also known that some of Tait's principles cannot be true in general. For instance the unknotting number cannot  always be determined by considering minimal projections. See \cite{Bleiler} and \cite{Nakanishi}. It was made clear in \cite{QW1} that Tait was mostly interested in -achirality, while Mary Haseman was the first to seriously address the +achirality question. 

\vskip.1in

We call the following statement {\bf Tait's Conjecture on alternating $-$achiral knots}: 

Let $K$ be an alternating $-$achiral knot. Then there exist a minimal projection $\Pi$ of $K$ in $S^2 \subset S^3$ and an involution $\varphi: S^3 \rightarrow S^3$ such that:
\\ 1) $\varphi$ reverses the orientation of $S^3$;
\\ 2) $\varphi(S^2) = S^2$;
\\ 3) $\varphi (\Pi) = \Pi$; 
\\ 4) $\varphi$ has two fixed points on $\Pi$ and hence reverses the orientation of $K$.

The reasons for which the paternity of this ``conjecture" may be  attributed to Tait are expounded in \cite{QW1}. The purpose of this paper is to prove this statement. 

Before moving towards the proof, we pause to point out that there is a significant difference between +achirality and -achirality. Here is why. All alternating achiral knots are hyperbolic. It results from several beautiful theorems in three dimensional topology that, if an hyperbolic knot $K$ is $\pm$achiral, there exists a diffeomorphism $\Phi : S^3 \rightarrow S^3$ as in Definition 1.1 which is periodic. Let us call the minimal period of such diffeomorphisms  $\Phi$ the {\bf order of $\pm$achirality of $K$.} Now, the order of -achirality of an hyperbolic knot is always equal to 2. Roughly speaking, Tait's Conjecture states that for -achirality, there exists a $\Phi$ of minimal order (i.e. of order 2) which can be seen on a minimal projection of $K$. On the other hand,  the order of +achrality can be any power of 2. When the order is equal to $2^a$ with $a \geq 2$ there exist alternating $+$achiral knots which have no minimal projection on which the symmetry is visible. Examples were provided by Dasbach-Hougardy for $a = 2$. See \cite{DH}. These facts are related to the isomorphism problem between the two checkerboard graphs, which is treated in Section 7.

The main tool to address chirality and more generally symmetry questions about alternating knots is provided by Tait Flyping Conjecture, proved by William Menasco and Morwen Thistlethwaite \cite{MT}.

\begin{theorem}{(Flyping Theorem)}
Let $\Pi_1$ and $\Pi_2$ be two minimal oriented projections in $S^2$ representing the same isotopy class of oriented alternating links in $S^3$. Then one can transform $\Pi_1$ into $\Pi_2$ by a finite sequence of flypes and orientation preserving diffeomorphisms of $S^2$. 
\end{theorem}

Let us recall that a {\bf flype} is the transformation of projections introduced by Tait and represented in Figure~1.

\begin{figure}[ht]    
   \centering
    \includegraphics[scale=0.45]{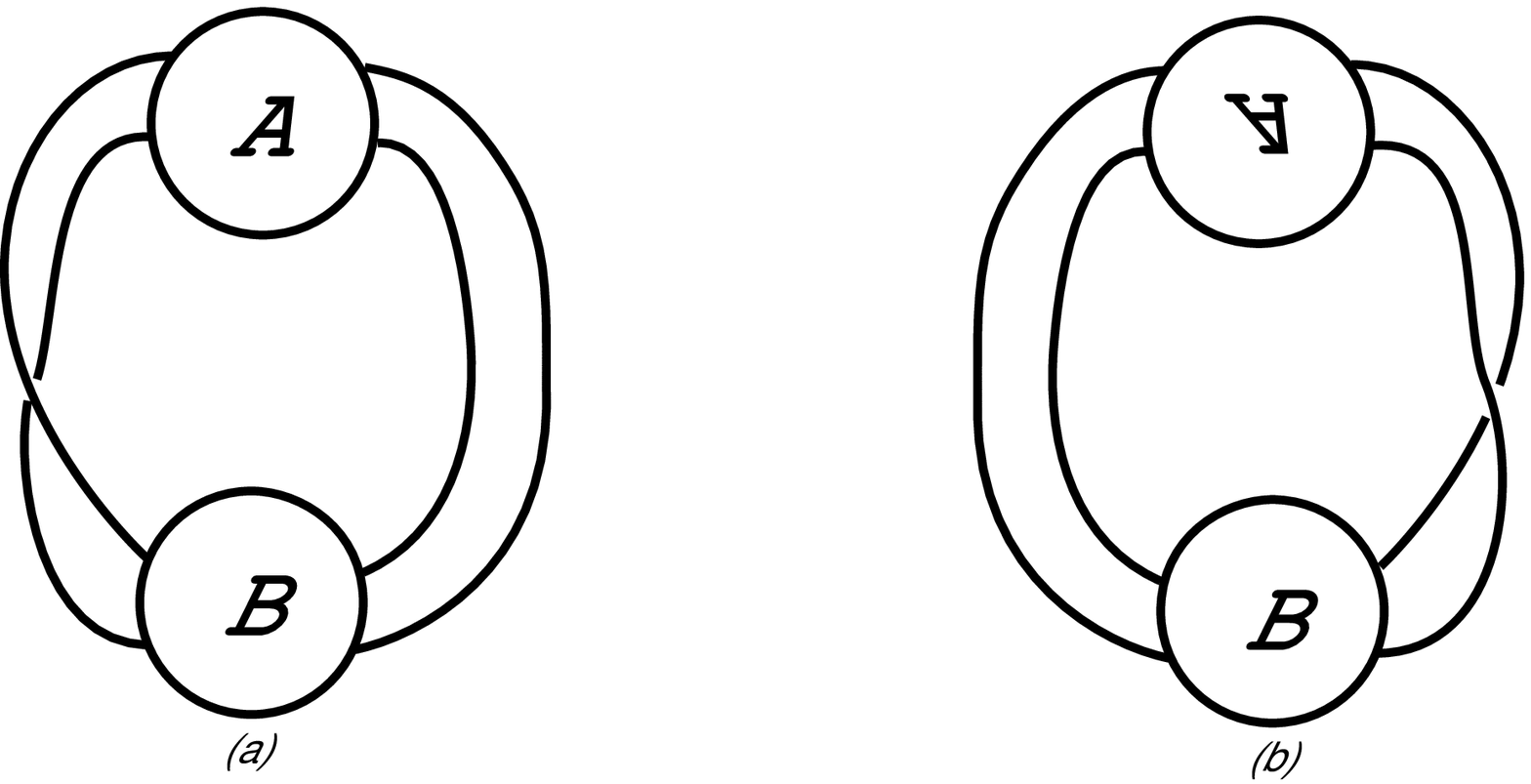}
\caption{A flype}
\end{figure}

{\bf Notation.} Let $\Pi$ be an oriented knot projection in $S^2$. We denote by $\widehat{\Pi}$ the mirror image of $\Pi$ through the sphere $S^2$ containing $\Pi$ and by $+\Pi$ resp. $-\Pi$ the projection obtained by preserving resp. reversing the orientation of $\Pi$.

The following result is a consequence of the Flyping Theorem.

\begin{theorem}{(Key Theorem)} Let $K$ be an alternating oriented knot. Let $\Pi$ be an oriented minimal projection of $K$. Then $K$ is $\pm$achiral if and only if one can proceed from $\Pi$ to $\pm\widehat{\Pi}$ by a finite sequence of flypes and orientation preserving diffeomorphisms of $S^2$.
\end{theorem}

Roughly speaking, to produce a proof of Tait's Conjecture we therefore have  to prove that in the case of $-$achirality we can always find a minimal projection for which no flypes are needed in the Key Theorem.

To achieve this goal, we have to know where flypes can be. The answer is provided by the paper \cite{QW2} in which, following Francis Bonahon and Laurent Siebenmann, the authors decompose in a canonical way  a knot projection $\Pi$ into {\bf jewels} and {\bf twisted band diagrams}. It turns  out that the flypes are all situated in the latter. Moreover Bonahon-Siebenmann's decomposition is (partly) coded by a finite tree. All minimal projections of a given alternating knot $K$ produce the same tree (up to canonical isomorphism). We call it the {\bf structure tree} of $K$ and denote it by ${\cal{A}} (K)$. The achirality of $K$ produces an automorphism $\varphi$ of ${\cal{A}} (K)$. The proof continues by determining the possible fixed points of $\varphi$ and by providing an argument in each case. 

{\bf The content of the paper is the following.} 

In Section 2 we recall the salient points about jewels and twisted band diagrams. Under mild assumptions (stated in four hypotheses) a (non-necessarily alternating) link projection $\Pi$  can be decomposed by Haseman circles (circles which cut $\Pi$ in four points) in a canonical way. This gives rise to the canonical decomposition of $\Pi$ in jewels and twisted band diagrams.
\\
In Section 3 we make explicit the position of flypes and we construct the structure tree $\mathcal {A} (K)$.
\\
In Section 4 we present the broad lines of the proof of Tait's Conjecture. Since $K$ is achiral,  the Key Theorem implies that there exists an automorphism $\varphi : \mathcal {A} (K) \rightarrow \mathcal {A} (K)$ which reflects the achirality of $K$. The proof then proceeds by an analysis of the possible fixed points of $\varphi$. They correspond to an invariant subset of $(S^2 ; \Pi)$ which can be a priori a jewel, a twisted band diagram or a Haseman circle. In the second part of Section 4 we prove that a twisted band diagram cannot be invariant.
\\
In Section 5 we prove Tait's Conjecture if a jewel is invariant. 
\\
In Section 6 we investigate achirality when a Haseman circle $\gamma$ is invariant. This circle separates $\Pi$ in two tangles. As they are exchanged by $\varphi$, they differ by their position with respect to $\gamma$ and by flypes. A detailed study of the various possibilities is thus undertaken. A proof of Tait's Conjecture follows, as well as some results about +achirality.
\\
In Section 7 we apply our methods to Kauffman's and Kauffman-Jablan's Conjectures about checkerboard graphs. We show that Kauffman's Conjecture is true for $-$achiral alternating knots and that Kauffman-Jablan's Conjecture is not true.

\section{The canonical decomposition of a projection}

In this section we do not assume that link projections are alternating.

\subsection{Diagrams}

\begin{definition}
A {\bf planar surface} $\Sigma$ is a compact connected surface embedded in the 2-sphere $S^2$. We denote by $v$ the number of connected components of the boundary $b\Sigma$ of $\Sigma$. 
\\
We consider compact graphs $\Gamma$ embedded in $\Sigma$ and satisfying the following four conditions:
\\ 1) vertices of $\Gamma$ have valency 1 or 4.
\\ 2) let $b\Gamma$ be the set of vertices of $\Gamma$ of valency 1. Then $\Gamma$ is properly embedded in $\Sigma$, i.e. $b\Sigma \cap \Gamma = b\Gamma$.
\\ 3) the number of vertices of $\Gamma$ contained in each connected component of $b\Sigma$ is equal to 4.
\\ 4) a vertex of $\Gamma$ of valency 4 is called a {\bf crossing point}. We require that at each crossing point an over and an under thread be chosen and pictured as usual. We denote by $c$ the number of crossing points.
\end{definition}

\begin{definition}
The pair $D = (\Sigma , \Gamma)$ is called a {\bf diagram}.
\end{definition}

\begin{definition}
 A {\bf singleton} is a diagram diffeomorphic to Figure 2(a).
\end{definition}
\begin{definition}
 A {\bf twisted band diagram} is a diagram diffeomorphic to Figure 2(b).
\end{definition}

\begin{figure}[ht]    
   \centering
        \includegraphics[scale=0.3]{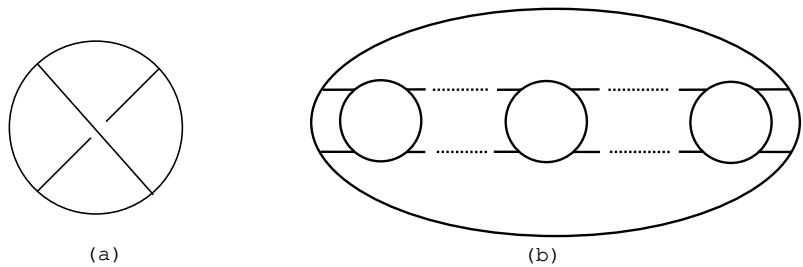}
\caption{(a) a singleton    (b) a twisted band diagram}
\end{figure}

\vskip.2in

The sign of a crossing point sitting  on a band is defined according to Figure 3.

\vskip.2in

\begin{figure}[ht]    
   \centering
    \includegraphics[scale=0.45]{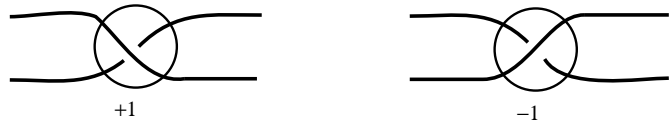}
\caption{The sign of a crossing on a band}
\end{figure}

{\bf First hypothesis.} Crossing points sitting side by side along the same band have the same sign. In other words we assume that Reidemeister Move 2 cannot be applied to reduce the number of crossing points along a band. 

Let us picture again a twisted band diagram as illustrated in Figure 4.

\begin{figure}[ht]    
   \centering
    \includegraphics[scale=0.35]{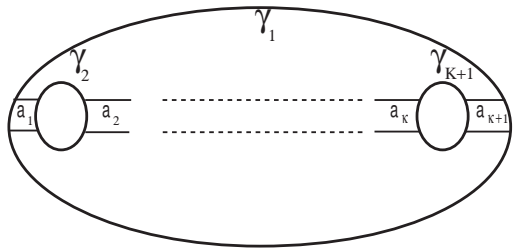}
\caption{A twisted band diagram with intermediate weights}
\end{figure}

In  Figure 4  the boundary components of $\Sigma$ are denoted by $\gamma_1 , \dots , \gamma_v$ where $v \geq 1$. The $a_i$ are integers. $\vert a_i \vert$ denotes the number of crossing points sitting side by side between $\gamma_{i-1}$ and $\gamma_i$. The sign of $a_i$ is the sign of the crossing points. The integer $a_i$ will be called an {\bf intermediate weight}. The corresponding portion of the diagram is called a {\bf twist}.

{\bf Second hypothesis.} If $v = 1$  we assume that $\vert a_1 \vert \geq 2$. If $v = 2$ we assume that $a_1$ and $a_2$ are not both 0. 

{\bf Remark.} Using flypes and then Reidemeister Move 2, we can reduce the number of crossing points of a twisted band diagram in such a way that either $a_i \geq 0$ for all $i = 1 , \dots , v$ or $a_i \leq 0$ for all $i = 1 , \dots , v$.

This reduction process is not quite canonical, but any two diagrams reduced in this manner are equivalent by flypes. This is enough for our purposes.

{\bf Third hypothesis.} We assume that in any twisted band diagram, all the non-zero $a_i$  have the same sign.

{\bf Notation.} The sum of the $a_i$ is called the {\bf weight} of the twisted band diagram and is denoted by $a$. If $v \geq 3$ we may have $a = 0$.

\subsection {Haseman circles}

\begin{definition}
 A {\bf Haseman circle} for a diagram $D = (\Sigma , \Gamma)$ is a circle $\gamma \subset \Sigma$ meeting $\Gamma$ transversally in four points, far from  crossing points. A Haseman circle is said to be {\bf compressible} if:
\\
i) $\gamma$ bounds a disc $\Delta$ in $\Sigma$.
\\
ii) There exists a properly embedded arc $\alpha \subset \Delta$ such that $\alpha \cap \Gamma = \emptyset$ and such that $\alpha$ is not boundary parallel. The arc $\alpha$ is called a {\bf compressing arc} for $\gamma$.
\end{definition}

{\bf Fourth hypothesis.} Haseman circles are  incompressible.

Two Haseman circles are said to be {\bf parallel} if they bound an annulus $A \subset \Sigma$ such that the pair $(A , A \cap \Gamma)$ is diffeomorphic to Figure 5.
\begin{figure}[ht]    
   \centering
    \includegraphics[scale=0.35]{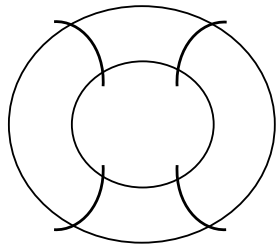}
\caption{Parallel Haseman circles}
\end{figure}

Analogously, we define a Haseman circle $\gamma$ to be {\bf boundary parallel} if there exists an annulus $A \subset \Sigma$ such that:
\\
1) the boundary $bA$ of $A$ is the disjoint union of $\gamma$ and a boundary component of $\Sigma$;
\\
2) $(A , A \cap \Gamma)$ is diffeomorphic to Figure 6.

\begin{definition}
 A {\bf jewel} is a diagram which satisfies the following four conditions:
\\ a) it is not a singleton.
\\ b) it is not a twisted band diagram with $v = 2$ and $a = \pm 1$.
\\ c) it is not a twisted band diagram with $v = 3$ and $a = 0$.
\\ d) every Haseman circle in $\Sigma$ is either boundary parallel or bounds a singleton.
\end{definition}

\begin{figure}[ht]    
   \centering
    \includegraphics[scale=1.]{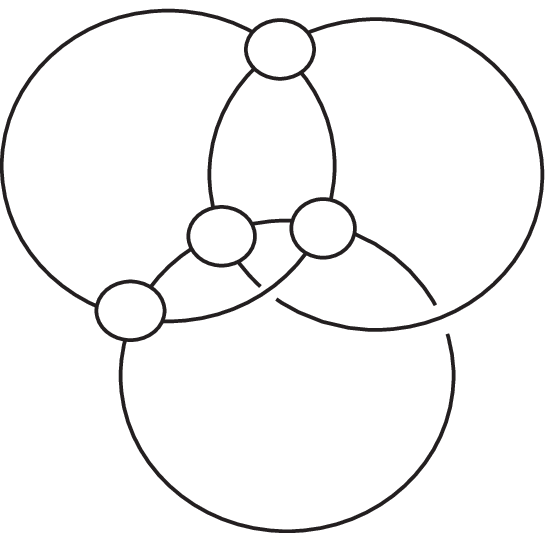}
\caption{A jewel}
\end{figure}

{\bf Comment.} The diagrams listed in a), b) and c) satisfy condition d) but we do not wish them to be jewels. As a consequence, a jewel is neither a singleton nor a twisted band diagram.

\subsection{Families of Haseman circles for a projection}

\begin{definition}
 A {\bf link projection} $\Pi$ (also called a projection for short) is a diagram in $\Sigma = S^2$. 
\end{definition}

{\bf Fifth hypothesis.} The projections we consider are connected and prime.

\begin{definition} 
 Let $\Pi$ be a link projection. A {\bf family of Haseman circles} for $\Pi$ is a set of Haseman circles satisfying the following conditions:
\\
1. any two circles are disjoint.
\\
2. no two circles are parallel.
\end{definition}

Note that a family is always finite, since a projection has a finite number of crossing points.

Let $\mathcal{H} = \{ \gamma_1 , ... , \gamma_n \}$ be a family of Haseman circles for $\Pi$. Let $R$ be the closure of a connected component of $S^2 \setminus  \bigcup_{i=1}^{i=n} \gamma_i$. We call the pair $(R , R \cap \Gamma)$ a {\bf diagram} of $\Pi$ determined by the family $\mathcal{H}$.

\begin{definition}
 A family $\mathcal{C}$ of Haseman circles is an {\bf admissible family} if each diagram determined by it is either a twisted band diagram or a jewel. An admissible family is {\bf minimal} if the deletion of any circle  transforms it into a family which is not admissible. 
\end{definition}

The next theorem is the main structure theorem about link projections proved in \cite{QW2}. It is essentially due to Bonahon and Siebenmann. 

\begin{theorem}{(Existence and uniqueness theorem of minimal admissible families)} Let $\Pi$ be a link projection in $S^2$. Then:
\\ i) there exist minimal admissible families for $\Pi$.
\\ ii) any two minimal admissible families are isotopic, by an isotopy which respects $\Pi$.
\end{theorem}

\begin{definition}
 ``The" minimal admissible family will be called the {\bf canonical Conway family} for $\Pi$ and denoted by $\mathcal{C}_{can}$. The decomposition of $\Pi$ into twisted band diagrams and jewels determined by $\mathcal{C}_{can}$ will be called the {\bf canonical decomposition } of $\Pi$.
\end{definition}

It may happen that $\mathcal{C}_{can}$ is empty. The next proposition tells us when this occurs.

\begin{proposition} Let $\Pi$ be a link projection. Then $\mathcal{C}_{can} = \emptyset$ if and only if $\Pi$ is either a jewel with empty boundary (i.e. $v = 0$) or the minimal projection of the torus knot/link of type $(2, m)$.
\end{proposition}

{\bf Comment.} A jewel with empty boundary is nothing else than a polyhedron in John Conway's sense which is indecomposable with respect to tangle sum. The minimal projection of the torus knot/link of type $(2, m)$ can be considered as a twisted band diagram with $v = 0$.

\section{The position of flypes and  the structure tree}

\subsection {The position of flypes}

From now on, we assume that the projection $\Pi$ we are considering is alternating and minimal.
The next theorem is an immediate consequence of the machinery presented in Section 2. Its importance for us is that we  now can locate with precision where the flypes happen. 
Thanks to the Flyping Theorem we can thus have a good idea about all minimal projections of an alternating link.

\begin{theorem}
(Position  of flypes.) Let $\Pi$ be a link projection in $S^2$ and suppose that a flype occurs in $\Pi$. Then, its active crossing point belongs to a twisted band diagram determined by $\mathcal{C}_{can}$. The flype moves the active crossing point
\\
1) either inside the twist to which it belongs,
\\
2) or  to another twist of the same band diagram.
\end{theorem}

{\bf Comments.} 1) The active crossing point is the one that moves during the flyping transformation. 
\\
2) Sometimes a flype of type 1) is called an {\bf inefficient flype} while one of type 2) is called {\bf efficient}. We are interested mainly in efficient flypes. 

\begin{definition}
We call the set of crossing points of the twists of a given twisted band diagram a {\bf flype orbit}.
\end{definition}

\begin{corollary}
A flype moves an active crossing point inside the flype orbit to which it belongs. Two distinct flype orbits are disjoint. 
\end{corollary}

The corollary can be interpreted as a loose kind of commutativity of flypes. Compare with \cite{Calvo}.

\subsection{The structure tree  ${\cal{A}} (K)$}

{\bf Construction of   ${\cal{A}} (K)$.}
Let $K$ be an alternating link and let $\Pi$ be a minimal projection of $K$. Let $\mathcal{C}_{can}$ be the canonical Conway family for $\Pi$. We construct the tree ${\cal{A}} (K)$ as follows. Its vertices are in bijection with the diagrams determined by $\mathcal{C}_{can}$. Its edges are in bijection with the Haseman circles of $\mathcal{C}_{can}$. The extremities of an edge (representing a Haseman circle $\gamma$)  are the vertices which represent the two diagrams which contain the circle $\gamma$ in their boundary. Since the diagrams are planar surfaces of a decomposition of the 2-sphere $S^2$ and since $S^2$ has genus zero, the graph we have constructed is a tree. This tree is ``abstract", i.e. it is not embedded in the plane. 

We label the vertices of  ${\cal{A}} (K)$ as follows. If a vertex represents a twisted band diagram we label it with the letter $B$ and by the weight $a$. If the vertex represents a jewel we label it with the  letter $J$. 

{\bf Remarks.} 1) The tree  ${\cal{A}} (K)$ is independent of the minimal projection chosen to represent $K$. This is an immediate consequence of the Flyping Theorem. Indeed, as we have seen,  the flypes  modify the decomposition of the weight $a$ of a twisted band diagram as the sum of intermediate weights, but the sum remains constant. A flype also modifies  the way in which diagrams are embedded in $S^2$. Since the tree is abstract a flype has no effect on it.  See \cite{QW2} Section 6. This is why we call it the {\bf structure tree} of $K$ (and not of $\Pi$). 
\\2) ${\cal{A}} (K)$   contains some information about the decomposition of $S^2$ in  diagrams determined by $\mathcal{C}_{can}$ but we cannot reconstruct the decomposition from it. However one can do better if no jewels are present. In this case the link (and its minimal projections) are called {\bf arborescent} by Bonahon-Siebenmann. They produce a planar tree which actually codes a given arborescent projection. See \cite{QW2} for details.
\\3) If $K$ is oriented, we do not  encode the orientation in ${\cal{A}} (K)$.

\begin{definition} (In the spirit of \cite{BS}.)
If all vertices of ${\cal{A}} (K)$ have a $B$ label the alternating knot $K$ is said to be {\bf arborescent}, while it is said to be {\bf polyhedral} if all its vertices have a $J$ label.  
\end{definition}

\section{First steps towards the proof of Tait's Conjecture}

\subsection{The automorphism of the structure tree}

Let $\widehat{K}$ be the mirror image of $K$. Let $\widehat{\Pi}$ be the mirror image of a minimal projection $\Pi$ of $K$. It differs from $\Pi$ by the sign of the crossings. Hence the tree ${\cal{A}} (\widehat{K})$ is obtained from ${\cal{A}} (K)$ by reversing the weight sign at each $B$-vertex. As abstract trees without signs at $B$-vertices, the two trees are canonically isomorphic (``equal").

{\bf Suppose now that $K$ is achiral (it does not matter whether + or $-$).}

The Key Theorem says that  there exists an isomorphism $\psi : (S^2 , K) \rightarrow (S^2 , \widehat{K})$ which is a composition of flypes and orientation preserving diffeomorphisms. This isomorphism induces an isomorphism  ${\cal{A}} (K) \rightarrow {\cal{A}} (\widehat{K})$. We  interpret it as an automorphism $\varphi : {\cal{A}} (K) \rightarrow {\cal{A}} (K)$ which, among other things, sends a $B$-vertex of weight $a$ to a $B$-vertex of weight $-a$. The Lefschetz Fixed Point Theorem implies that $\varphi : {\cal{A}} (K) \rightarrow {\cal{A}} (K)$ has fixed points.

{\bf Remark.} During the proof we shall see that $\varphi$ has just one fixed point.

The proof of Tait's Conjecture is divided in three cases, depending (a priori)  whether:
\\
1) a twisted band diagram is invariant (in fact we shall prove immediately below that this is impossible);  
\\
2) a jewel is invariant;
\\
3) a Haseman circle is invariant. This last case is  the most involved.

\subsection{A twisted band diagram cannot be invariant}

\begin{lemma}
Let $\Pi$ be an alternating projection and let $\gamma$ be a Haseman circle for $\Pi$. Choose one side of $\gamma$. Consider the four threads of $\Pi$ which cut $\gamma$. Label each thread with a + or a $-$ if at the next crossing  on the chosen side the thread will pass above or below. Then opposite threads along $\gamma$ have the same label.
\end{lemma}

{\bf Proof.} Consider the chessboard associated to $\Pi$. Then $\Pi$ is alternating according to the following rule. Consider a thread  slightly  before a crossing and move along the thread towards the crossing. If a given colour (say black) is on the  right during the move then the thread will be an overpass. Since opposite regions along $\gamma$ have the same colour  this proves the lemma.

{\bf Consequence of the lemma.} Let $\Pi$ be a minimal alternating projection. Consider its canonical Conway family $\mathcal{C}_{can}$. Let $D$ be a diagram (either a twisted band diagram or a jewel) determined by $\mathcal{C}_{can}$. Let $\gamma$ be a boundary component of $D$. Then it is always possible to attach a singleton on $\gamma$ in order to obtain an alternating projection $D^*$. This procedure is unique. The next lemma follows immediately.

\begin{lemma}
Suppose that the diagram $D$ considered above is invariant by $\psi$. Then the restriction of $\psi$ to $D$ extends to $D^*$. 
\end{lemma}

\begin{proposition}
A twisted band diagram cannot be invariant by $\psi$.
\end{proposition}

{\bf Proof.} Let $D$ be a twisted band diagram invariant by $\psi$. Let $\psi^* : D^* \rightarrow D^*$ be its extension to $D^*$. Since $D^*$ is the minimal alternating projection of a torus knot/link of type $( 2, m)$ which is chiral, the existence of such a diffeomorphism $\psi^*$ is impossible.

\begin{corollary}
Suppose that the edge of ${\cal{A}} (K)$ representing a Haseman circle $\gamma$ is invariant by $\varphi$. Then the extremities of this edge cannot be labeled one by a $B$ and the other by a $J$. 
\end{corollary}

{\bf Intermediate result in our quest to understand the achirality of alternating knots}. We are thus left with three possibilities:
\\
A) Existence of an invariant jewel: a jewel is invariant by $\varphi : {\cal{A}} (K) \rightarrow {\cal{A}} (K)$;
\\
B) Existence of a polyhedral invariant circle: a Haseman circle $\gamma$ is invariant by $\psi$ and both diagrams adjacent to $\gamma$ are jewels; we shall see in next section  that in this case the two diagrams adjacent to $\gamma$ are exchanged by $\psi$;
\\
C) Existence of an arborescent invariant circle: a Haseman circle $\gamma$ is invariant by $\psi$ and both diagrams adjacent to $\gamma$ are twisted band diagrams. By proposition 4.1 the two adjacent diagrams must be exchanged by $\psi$. 

{\bf Remark.} The terms we use (polyhedral or arborescent) do not mean that the whole projection is so, but only the region of the projection near the fixed point.

\section{Proof of Tait's Conjecture if a jewel is invariant}

\begin{proposition}
Let $\Pi$ be a projection  and    let $f : (S^2 , \Pi) \rightarrow (S^2 , \Pi)$ be a diffeomorphism. Then $f$ is isotopic (by an isotopy respecting $\Pi$) to a diffeomorphism of finite order. 
\end{proposition}

{\bf Proof.} Consider the cell decomposition of $S^2$induced by $\Pi$ and argue inductively on the dimension of cells.

\begin{lemma}(Tutte's Lemma)
Let $\Gamma \subset S^2$ be a finite graph embedded in the 2-sphere. Let $g : (S^2 , \Gamma) \rightarrow (S^2 , \Gamma)$ be an orientation preserving  diffeomorphism of finite order 
 which is the identity on an edge  of $\Gamma$. Then $g$ is the identity.
\end{lemma}

For a proof of the lemma see \cite{Tutte}.

\begin{proposition}
Let $K$ be an alternating $-$achiral knot. Let $\Pi$ be a minimal projection such that no flype is needed to transform $\Pi$ into $-\bar{\Pi}$. Then the projection $\Pi$ satisfies Tait's Conjecture.
\end{proposition}

{\bf Proof of the proposition.} If no flype is needed there is an orientation preserving diffeomorphism $f: S^2 \rightarrow S^2$ which transforms  $\Pi$ into $-\bar{\Pi}$. By the preceeding proposition we may assume that f is of  finite order. Since $\Pi$ is the image of a generic immersion of the circle $S^1$ in $S^2$,   the restriction of $f$  to $\Pi$ pulls back to a diffeomorphism $\hat f$ of finite order of $S^1$. Since $f$ reverses the orientation of $\Pi$, by Smith theory  $\hat f$ has 2 fixed points and is of order 2. Then, by Tutte's Lemma $f : S^2 \rightarrow S^2$ is of order 2.

{\bf The Filling Construction.} Let $\Pi$ be a knot projection and let $D$ be a jewel subdiagram of $\Pi$. Let $\gamma_1 , \dots , \gamma_k$ be the Haseman circles which are the boundary components of $D$. Each $\gamma_i$ bounds in $S^2$ a disc $\Delta_i$ which does not meet the interior of $D$.  The projection $\Pi$ cuts each $\gamma_i$ in 4 points. Inside $\Delta_i$ the projection $\Pi$ joins either opposite points or adjacent points. In the first case, we replace $\Delta_i \cap \Pi$ by a singleton. In the second case we replace $\Delta_i \cap \Pi$ as follows:

\begin{figure}[ht]    
   \centering
    \includegraphics[scale=0.4]{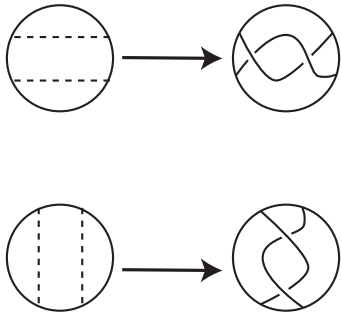}
\caption{Replacement}
\end{figure}

We obtain in this manner a new projection $\Pi^*$. The Filling Construction is such that $\Pi^*$ is again the projection of a knot. Using Lemma 4.1 we can choose the over/under crossings in $\Delta_i \cap \Pi^*$ in such a way that $\Pi^*$ is alternating. An important fact is that $\Pi^*$ has no place for a flype. 
\\
Suppose now that the knot $K$ represented by $\Pi$ is + or -achiral. Let $\psi : \Pi \rightarrow \widehat {\Pi}$ be as before a composition of flypes and diffeomorphisms produced by the Key Theorem and suppose that $D$ is invariant by $\psi$. We claim that  $\Pi^*$ is the projection of a + or -achiral knot as $\Pi$ was and that the restriction of $\psi$ to $D$ extends to 
$\psi^* : \Pi^* \rightarrow \widehat {\Pi}^*$.
Since $\Pi^*$ has no place for flypes, $\psi^*$ is  a diffeomorphism. 
\\
We also see  that $\psi^*$ (and hence also $\psi$) leaves no boundary component of $D$ invariant. This proves that in Case B) above the two jewels adjacent to $\gamma$ must be exchanged.

\begin{proposition}
Let $K$ be an alternating $-$achiral knot. Let $\Pi$ be a minimal projection of $K$ such that a jewel diagram is invariant by $\psi$. Then Tait's Conjecture is true for $\Pi$. 
\end{proposition}

{\bf Proof of the proposition.} Let $D$ be the jewel invariant by $\psi$. Let ${\Pi}^*$ be the projection obtained by the Filling Construction   above and let ${\psi}^* : {\Pi}^*  \rightarrow {\Pi}^*$ be the automorphism induced by $\psi$. It results from the construction of $\psi^*$ and from Proposition 5.2 that $\psi^*$ is an involution that leaves no Haseman circle invariant. Let $k = 2l$ be the number of these circles. We number the circles $\gamma_1 , \dots , \gamma_{2l}$ in such a way that $\psi^*(\gamma_i) = \gamma_{2l+1-i}$ for $i = 1 , \dots , 2l$. Let $D_i = (\Delta_i , \Gamma_i)$ be the diagram contained in the disc $\Delta_i$ bounded by $\gamma_i$. The diffeomorphism $\psi$ sends $\Delta_i$ onto $\Delta_{2l+1-i}$. It follows from the Key Theorem that $\psi (D_i)$ is flype equivalent to $D_{2l+1-i}$. Hence let us perform flypes in $D_{2l+1-i}$ to obtain $\hat {D}_{2l+1-i}$ in such a way that $\psi (D_i) = \hat {D}_{2l+1-i}$. Since $\psi^*$ is an involution we have $\psi(\hat {D}_{2l+1-i}) = D_i$. Thus we have constucted a projection which is invariant by an involution. 

{\bf Remarks on the proof.} 1) We have modified the initial projection by equivariant  flypes in the discs $\Delta_i$. This was possible since the involution $\psi^*$ acts freely on the family of  discs $\Delta_i$. 
\\
2) It is crucial in the proof that $\Pi$ is the projection of a knot (not of a link with several components). It is also crucial that we study -achirality (i.e. that the orientation of the knot is reversed).

{\bf Comments.} 1) Proceeding backwards, we see that projections satisfying Tait's Conjecture with an invariant jewel $D$ are obtained as follows. Start with a $-$achiral projection $D^*$ which is a jewel. Consider the automorphism $\psi^* : D^* \rightarrow D^*$ which realises the $-$achirality symmetry. Replace then the singletons by arbitrary diagrams in an equivariant way. This construction is well known to specialists and was certainly known to John Conway when he wrote his celebrated paper \cite{Conway}.
\\
2) The proof of Tait's Conjecture in Case B) (existence of an invariant polyhedral circle) follows the same lines as the proof in Case A). We just have to replace the invariant jewel $D$ by the union of the two jewels which are adjacent to the invariant Haseman circle $\gamma$. 
\\
3) We shall study in a future publication \cite{EQ} the status of +achiral alternating knots when a jewel is invariant. Among other things, we shall see that the order of achirality can be any power of $2$ and that this symmetry is not always realised by a diffeomorphism on a minimal projection (in other words, flypes may be needed). The main reason is that, when the order of achirality is not prime (i.e. equal to $2^a$, with $a \geq 2$), the permutation induced by $\psi^*$ on the discs $\Delta_i$ may have cycles of different lengths. Of course, this cannot happen if the order of achirality is equal to $2$; this is one reason why Tait's Conjecture is true for -achirality.

\section{Achirality when a Haseman circle is invariant}

{\bf In  this section we suppose that $K$ is $\pm$achiral}. We shall return to $-$achirality when needed. Our study applies to invariant Haseman circles, either polyhedral or arborescent.

\begin{figure}[ht]    
   \centering
    \includegraphics[scale=0.3]{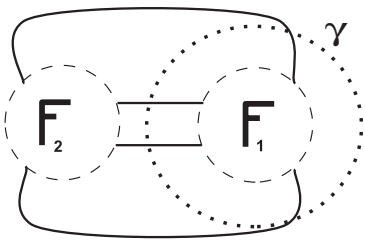}
\caption{The partition of $\Pi$ induced by the splitting of $\gamma$}
\end{figure}

\subsection{The partition of $\Pi$ induced by the invariant Haseman circle}

To prepare a further study of +achirality, we investigate the following situation, which includes Cases B) and C).  The origin of this study goes back to Mary Haseman.  
\\
We have a Haseman circle $\gamma$ invariant by $\psi$. The two  diagrams $D_1$ and $D_2$ adjacent to $\gamma$ are both twisted band diagrams or both jewels; they are exchanged by $\psi$. To better see what happens, we choose to picture the situation as shown in Figure 8, with the circle $\gamma$ split into two parallel circles. This produces a partition of $\Pi$ in two tangles. 

\begin{theorem}
Let $K$ be an oriented $\pm$achiral knot.  Suppose that $K$ has a minimal projection with an invariant Haseman circle $\gamma$. Then, up to a global change of orientation, $K$ admits a minimal projection of Type I or II, as shown in Figures 9 and 10. 
\end{theorem}

\begin{figure}[h]
\begin{minipage}[b]{0.5\linewidth}
  \centering
  \includegraphics[scale=0.3]{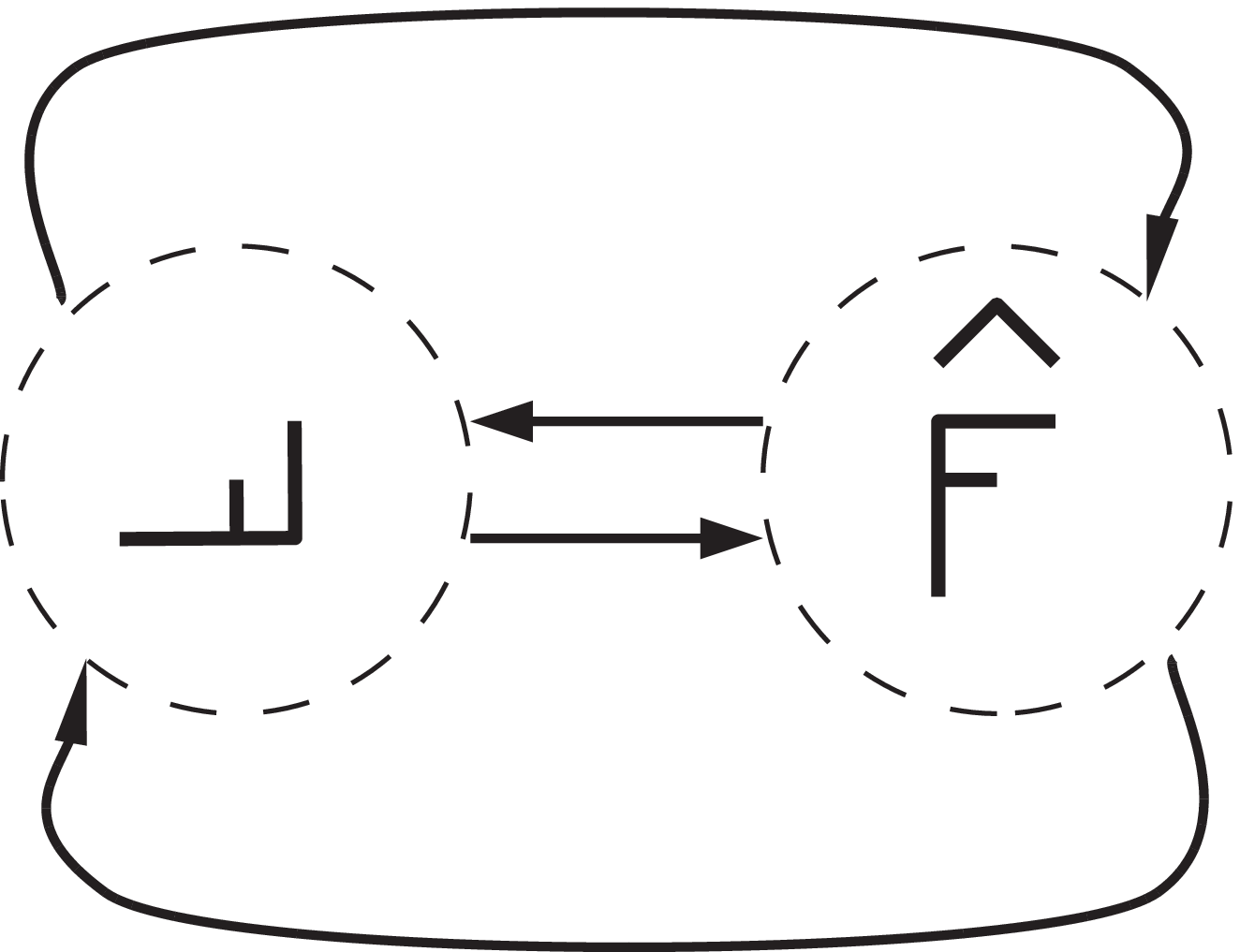}
  \caption{Projection of Type I}
  \end{minipage}
\hspace{0.5cm}
\begin{minipage}[b]{0.5\linewidth}
  \centering
  \includegraphics[scale=0.3]{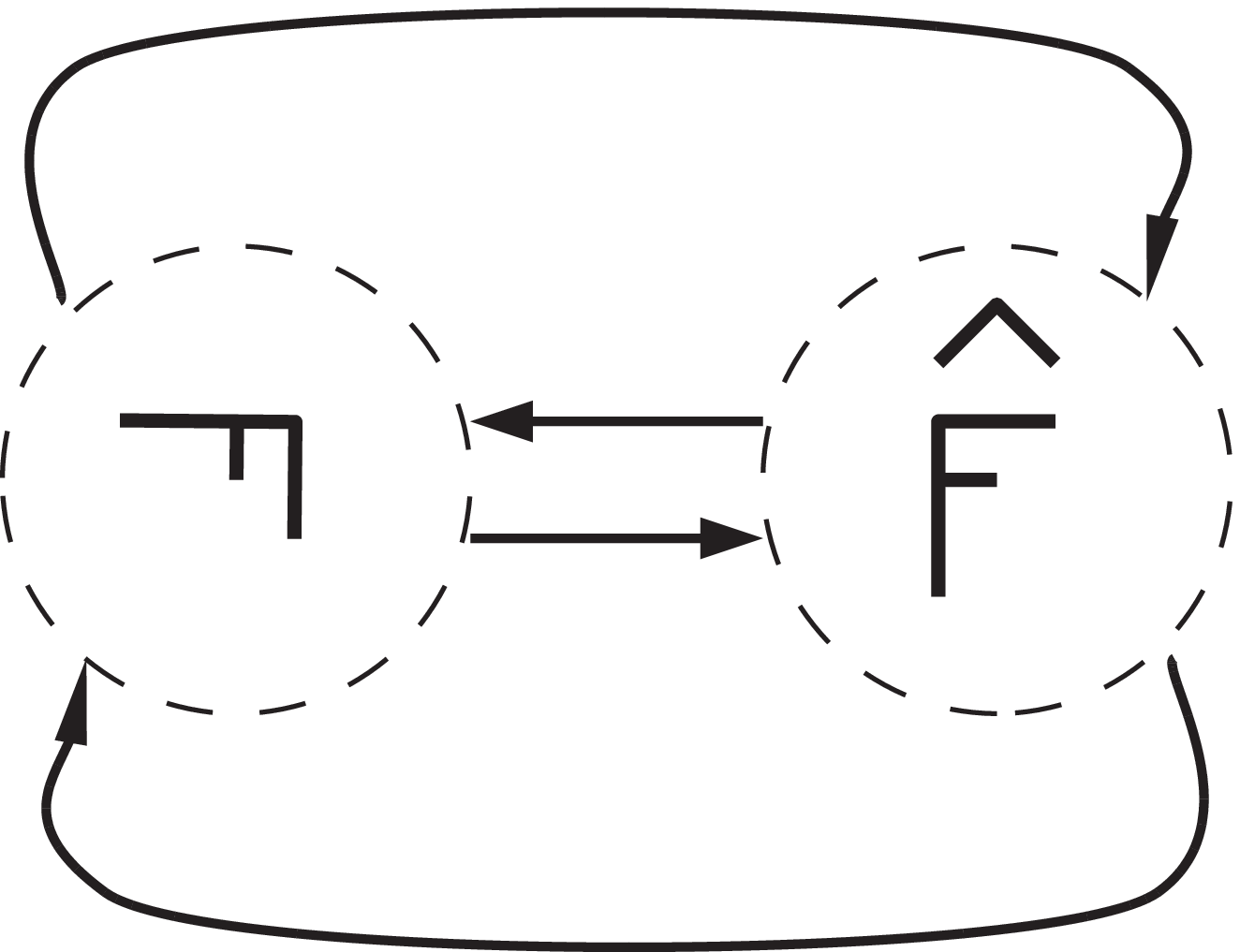}
    \caption{Projection of Type II}
    \end{minipage}
\end{figure}

{\bf Explanations about the figures}
\\
i) The symbol $F$ denotes the tangle inside the left-hand side circle;
\\
ii) the symbol $\widehat{F}$ denotes the image of the tangle $F$ by the reflection through the projection plane;
\\
iii) in the Type II picture, the right-hand tangle is obtained from the left-hand one by a rotation of angle $\pi \over 2$ with an axis orthogonal to the projection plane, followed by a reflection through the projection plane;
\\
iv) in the Type I picture, the right-hand tangle is obtained from the left-hand one by the following three moves:
\\
1) a rotation of angle $\pi \over 2$ with an axis orthogonal to the projection plane;
\\
2) a rotation of angle $\pi$ with an axis in the projection plane;
\\
3) a reflection through the projection plane.

{\bf The theorem applies to invariant Haseman circles which are polyhedral as well as arborescent.}

{\bf Proof of the theorem.}

We shall advance in the proof step by step.

{\bf Step 1.}  We know that the achirality of $K$ produces an isomorphism $\psi : \Pi \rightarrow \widehat {\Pi}$, via the Key Theorem. In turn, $\psi$ induces an automorphism $\varphi : {\cal{A}} (K) \rightarrow {\cal{A}} (K)$ which, by hypothesis, has a fixed point in the middle of an edge $\epsilon$. This reflects the fact that $\psi$ leaves invariant the Haseman circle $\gamma$ corresponding to $\epsilon$. We know that $\psi$ has to exchange the two tangles $F_1$ and $F_2$ (see Figure 8) which have $\gamma$ for boundary. Let us denote one of them by $F$. 
\\
Now comes a rather subtle point. The minimal projection $\Pi$ is not unique as it is acted upon by flypes. However:
\\
(i) the Haseman circle $\gamma$ is insensitive to flypes (because of  the theorem which determines the position of flypes); 
\\
(ii) flypes act independently in the two tangles determined by $\gamma$. 
\\
We summarise a consequence of these last arguments as follows.

{\bf Important observation.} If an alternating knot $K$ possesses a minimal projection of Type I or II, then this projection is unique up to flypes occuring in $F$ and/or in $\widehat {F}$.

Let $\sigma : (S^2 , \Pi) \rightarrow (S^2 , \widehat {\Pi)}$ be a diffeomorphism which exchanges the two tangles and  is part of $\psi$ (remember that $\psi$ is a composition of flypes and diffeomorphisms).  Since $K$ is achiral of some sort, there must exist by the Key Theorem flypes which transform $\sigma (F)$ into $\widehat {F}$. The point is that, without knowing exactly what $\sigma$ can be, the position of $\widehat {F}$ with respect to $F$ is not arbitrary. The theorem  says that the only possibilities are shown in projections of Type I and II. 
\\
The detailed study of the potential diffeomorphisms $\sigma$ will be undertaken in  later subsections. 

{\bf Step 2.} Neglecting orientations there are a priori eight possible partitions of $\Pi$. Indeed, keeping the tangle $\widehat {F}$ fixed, there are eight ways to place $F$ in the left-hand side tangle. Notice that the dihedral group $D_4$ acts freely and transitively on these left-hand side tangles, as the group of symmetries of a circle with four marked points.

\begin{figure}[ht]    
   \centering
    \includegraphics[scale=0.3]{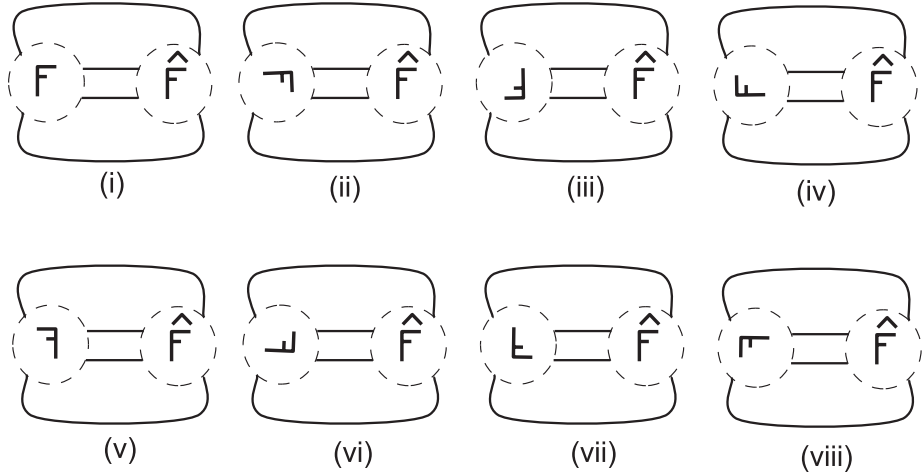}
\caption{The eight possible configurations}
\end{figure}

{\bf Step 3.} The cases (i), (iii), (v), (vii) are impossible. Proof:
\\
The boundary circle of a tangle is cut transversely in four points by the projection $\Pi$. Traditionally, these four points are located in North-West (NW), North-East (NE), South-East (SE) and South-West (SW). In each tangle, the projection $\Pi$ connects  these four points in pairs. We call this connection a ``connection path". A priori, there are three possible such paths:
\\
NW is connected to NE, and SW is connected to SE. We denote this  connection path by $H$, for ``horizontal".
\\
NW is connected to SW and NE is connected to SE. We denote this connection path by $V$, for ``vertical".
\\
NW is connected to SE and SW is connected to NE. We denote this connection path by $X$, for ``crossing".\\

\begin{figure}[ht]    
   \centering
    \includegraphics[scale=0.4]{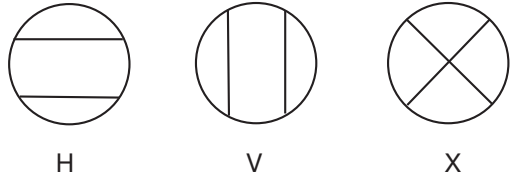}
\caption{Description of the connection paths H, V and X}
\end{figure}

Now, in the figures (i), (iii), (v) and (vii) the positions of $F$ and $\widehat {F}$ differ essentially by a rotation of angle $\pi$. This leaves the connection paths invariants. In other words, the connection paths are the same in the two tangles. A quickly drawn picture shows that, in  each  of the three cases (two $H$'s, or two $V$'s, or two $X$'s), the projection $\Pi$ represents a link with at least two components and not a knot.

{\bf Step 4.} Figure (ii) is isomorphic to (iv) and (vi) is isomorphic to (viii). Indeed, we can transform (ii) into (iv) and (vi) into (viii) by a rotation of angle $\pi$, with axis in the projection plane as indicated in Figure 13. The rotation axis is represented by a dotted line.\\

\begin{figure}[ht]    
   \centering
    \includegraphics[scale=0.35]{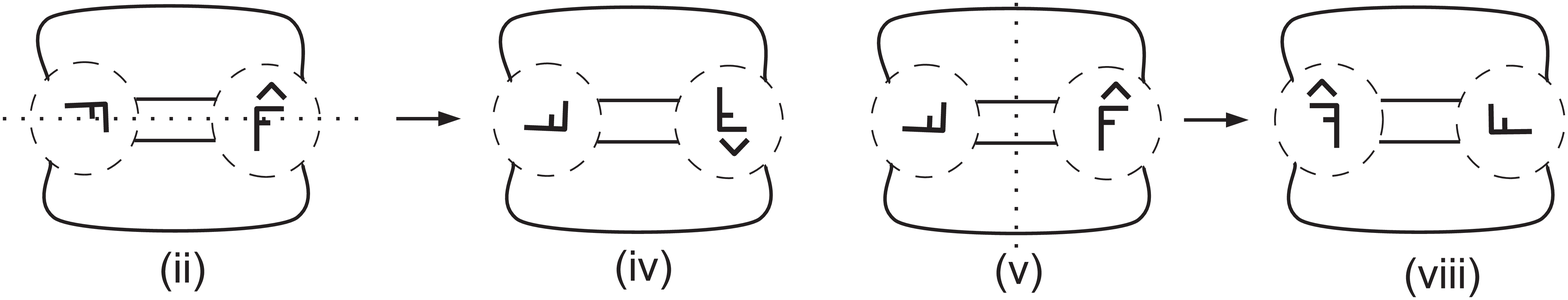}
\caption{}
\end{figure}
 
{\bf Step 5.} Orientations of the strands. We are thus left with the two figures (ii) and (vi). We wish to orient the four strands which connect the two tangles. We number them from 1 to 4, beginning with the top one. We say that an orientation of the stands, represented by an arrow on each strand, is a parallel one, if there exist two consecutive strands which have arrows pointing in the same direction. We claim that a parallel orientation of the strands is impossible. To prove this, we consider the connection paths inside the tangles. We observe that the two connection paths in $F$ and in $\widehat {F}$ differ by a rotation of angle $\pi \over 2$. If one is $X$, then so is the other. We have already seen that this is impossible. If one is $H$ then the other is $V$. But this is incompatible with a parallel orientation. 
\\
Hence, up to global change of orientations, there is only one possible orientation of the strands. It is the one in which two consecutive arrows have the same direction. In this case, one connection path is an $H$ while the other  is a $V$.  

{\bf End of proof of the theorem.}

Note that the theorem does not answer the achirality question of projections of Type I or II. This will be done in the next subsections.

\subsection{Achrality of projections of Type I}

\begin{definition}
Let $F_1$ and $F_2$ be two tangles. We say that they are {\bf flype equivalent} if one can obtain $F_2$ from $F_1$ by a sequence of flypes leaving the boundary circle fixed. We write $F_1 \sim F_2$ for this relation.
\end{definition}

Let $F$ be a tangle. We denote by:
\\
$F^*$ the tangle obtained from $F$ by a rotation of angle $\pi$, with an axis orthogonal to the projection plane.
\\
$F^h$ the tangle obtained from $F$ by a rotation of angle $\pi$, with an horizontal axis contained in the projection plane.
\\
$F^v$ the tangle obtained from $F$ by a rotation of angle $\pi$, with a vertical axis contained in the projection plane. 
\\
See Figure 14.\\

\begin{figure}[ht]    
   \centering
    \includegraphics[scale=0.65]{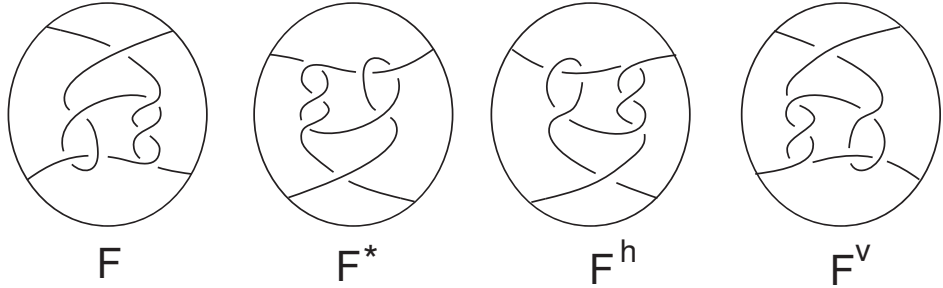}
\caption{A picture of the tangles $F$, $F^*$, $F^h$ and $F^v$}
\end{figure}

\begin{proposition}
Let $K$ be an alternating knot which possesses a minimal projection of Type I. Then:
\\
(i) $K$ is +achiral if and only if $F \sim F^*$;
\\
(ii) $K$ is -achiral if and only if $F \sim F^h$ or $F \sim F^v$.
\\
In particular, if $K$ is -achiral, then $K$ possesses also a projection of Type II. 
\end{proposition}

{\bf Proof of the proposition.}

The argument will be as follows. 
\\
We know that $\psi$ transforms $\Pi$ into $\widehat{\Pi}$ and leaves $\gamma$ invariant. Besides flypes, $\psi$ contains a   diffeomorphism $S^2 \rightarrow S^2$ which leaves $\gamma$ invariant and exchanges the two hemispheres. Now look at Figure 8, which contains the Haseman circle split in two, thus producing two tangles. Let us call this structure the {\bf frame} of a projection of Type I or II. The content of the circles (i.e  the tangles themselves) is not part of the structure. The frame consists of the two circles, connected by four strands. We denote it by $\Omega$ and consider it as embedded in $S^2$. 
\\
We are looking for automorphisms $\sigma : (S^2 ; \Omega) \rightarrow (S^2 ; \Omega)$ of the frame which exchange the two circles. We produce a complete list of such automorphisms (up to isotopy). Then we search among them for those which could occur in the symmetry $\psi$. These eligible ones set conditions on the tangle $F$.

We represent $S^2$ as an terrestrial  globe. We single out the Equator $E$ and we mark on $E$ the  degrees of longitude from $0^{\circ}$ to $360^{\circ}$. We trace two great circles $C_1$ and $C_2$ on $S^2$. Both go through the North Pole $P$ and the South Pole $S$  and are orthogonal to $E$.  The great circle $C_1$ intersects $E$ at the longitudes  $0^{\circ}$ and $180^{\circ}$ while $C_2$ intersects $E$ at $90^{\circ}$ and $270^{\circ}$. 

\begin{figure}[ht]    
   \centering
    \includegraphics[scale=0.4]{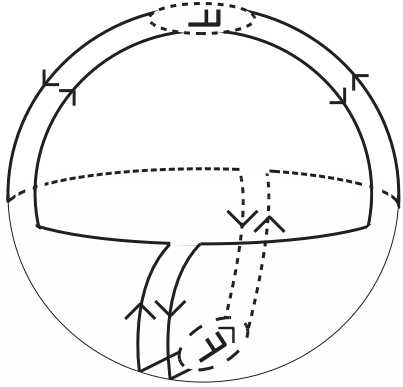}
\caption{A picture of $\Pi$ on the 2-sphere}
\end{figure}

By hypothesis, the invariant Haseman circle $\gamma$ divides $\Pi$ into two tangles $F$ and $\widehat {F}$.   The Equator $E$ is the Haseman circle $\gamma$. The tangle $F$ is contained in the Northern Hemisphere $NH$ and $\widehat {F}$ is contained in the Southern Hemisphere $SH$.  The ``band" which goes from $\gamma$  to $F$ is essentially concentrated along $C_1 \cap NH$ and the ``band" which goes from $\gamma$ to $\widehat {F}$  is concentrated along $C_2 \cap SH$. The tangle $F$ itself is near the North Pole, while the tangle $\widehat{F}$ is near the South Pole. See Figure 15.

We consider five rotation axes. They are among the symmetry axes of the octahedron drawn on $S^2$ consisting of $E$, $C_1$ and $C_2$. A priori, they are the only admissible ones for the picture we consider. 

The lines $L_1$, $L_2$, $L_3$ and $L_4$ are contained in the equatorial plane in such a way that:
\\
 $L_1$  cuts the Equator at $135^{\circ}$ and $315^{\circ}$.
\\
 $L_2$  cuts the Equator at $45^{\circ}$ and $225^{\circ}$.
\\
$L_3$  cuts the Equator at $0^{\circ}$ and $180^{\circ}$.
\\
$L_4$  cuts the Equator at $90^{\circ}$ and $270^{\circ}$.

$L_5$ is the vertical line which goes through the poles.

For $i = 1, \dots , 5$ we denote by $r^{L_i}_n$ the rotation of angle $2\pi \over n$, with  the line $L_i$ as axis. We denote by $\rho$ the reflection through the equatorial plane. 

It is easily checked that the only orthogonal transformations  $\sigma : (S^2 ; \Omega) \rightarrow (S^2 , \Omega)$ are the following eight ones:

(i) id (the identity)

(ii) $r^{L_1}_2$

(iii) $r^{L_2}_2$

(iv) $\rho \circ r^{L_5}_4 = R^{L_5}_4$

(v) $\rho \circ r^{L_5}_{-4} = R^{L_5}_{-4}$

(vi) $\rho \circ r^{L_3}_2 = R^{L_3}_2$

(vii) $\rho \circ r^{L_4}_2 = R^{L_4}_2$

(viii) $r^{L_5}_2$

We now visualise and analyse  each  of the eight symmetries. Since we are interested in achirality, we use the following principle.
\\
Let $\sigma$ be any of the eight symmetries. Then $\sigma$ is relevant for the achirality of the knot $K$ represented by the projection $\Pi$ if either:
\\
1) $\sigma$ preserves the orientation of $S^3$ (represented in the figures by ${\bf R}^3$) and sends $\Pi$ to $\widehat {\Pi}$;
\\
or  2) $\sigma$ reverses the orientation of $S^3$ and sends $\Pi$ to $\Pi$.
\\
But if $\sigma$ reverses the orientation of $S^3$ and sends $\Pi$ to $\widehat {\Pi}$, then $\sigma$ is not pertinent.

{\bf Analysis  of (i), the identity symmetry}. This means that we can modify the tangle $F$ to the tangle $\widehat {F}$ by flypes only, leaving the boundary circle fixed. This is impossible. Here is why.
\\
If no flypes are possible in $F$, there is nothing to prove since $F$ contains crossings and since $\widehat {F}$ differs from $F$ by the over/under passes at crossings.
\\
Hence suppose by contradiction that flypes occur. We have seen in Section  3 that flypes  occur in twisted band diagrams. If a flype occurs in such a diagram $D$, the weight of $D$ is $\neq 0$. We know that weights are invariant by flypes, and that the weights of $\widehat {F}$ are the opposite of those of $F$. If the weight of $D$ is different from $0$, we reach a contradiction as $D$ does not move since $\gamma$ is fixed.  Otherwise, we can find a chain of enclosed twisted band diagrams $D = D_1 \supset D_2 \supset \dots \supset D_k$ such that the weight of $D_i$ is $0$ for $1 \leq i \leq k-1$ and the weight of $D_k$ is different from $0$. We reach a contradiction as above.

{\bf Analysis of (ii) $r^{L_1}_2$}. The action of $r^{L_1}_2$ is represented in Figure 16.

\begin{figure}[ht]    
   \centering
    \includegraphics[scale=0.35]{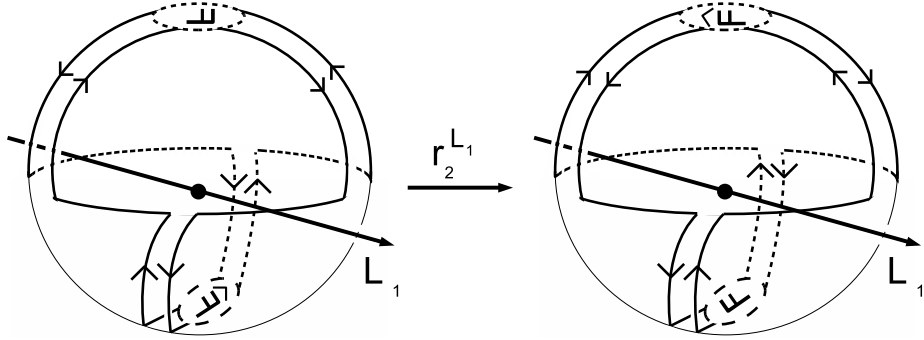}
\caption{The action of $r^{L_1}_2$}
\end{figure}

We see that we can transform $\Pi$ into $\widehat {\Pi}$ (actually into $-\widehat {\Pi}$) with flypes and the rotation $r^{L_1}_2$ if and only if $F \sim F^h$. Remark that $F \sim F^h$ implies that $\Pi$ is a Type II projection.

{\bf Analysis of (iii) $r^{L_2}_2$}. The action of  $r^{L_2}_2$ is represented in Figure 17. 

\begin{figure}[ht]    
   \centering
    \includegraphics[scale=0.35]{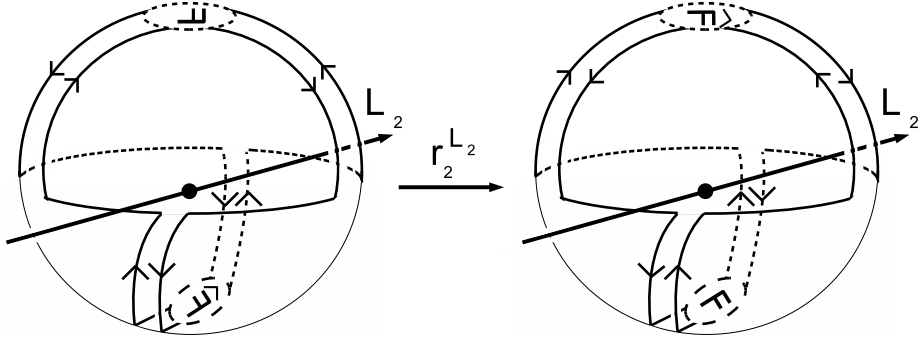}
\caption{The action of $r^{L_2}_2$}
\end{figure}

The action of this rotation is very similar to the preceding one. We see that we can transform $\Pi$ into $\widehat {\Pi}$ (actually into $-\widehat {\Pi}$) with flypes and the rotation  $r^{L_2}_2$ if and only if $F \sim F^v$. Remark that $F \sim F^v$ implies that $\Pi$ is a Type II projection.

{\bf Analysis of (iv) $R^{L_5}_4$}. The action of $R^{L_5}_4$ is represented in Figure 18.

\begin{figure}[ht]    
   \centering
    \includegraphics[scale=0.35]{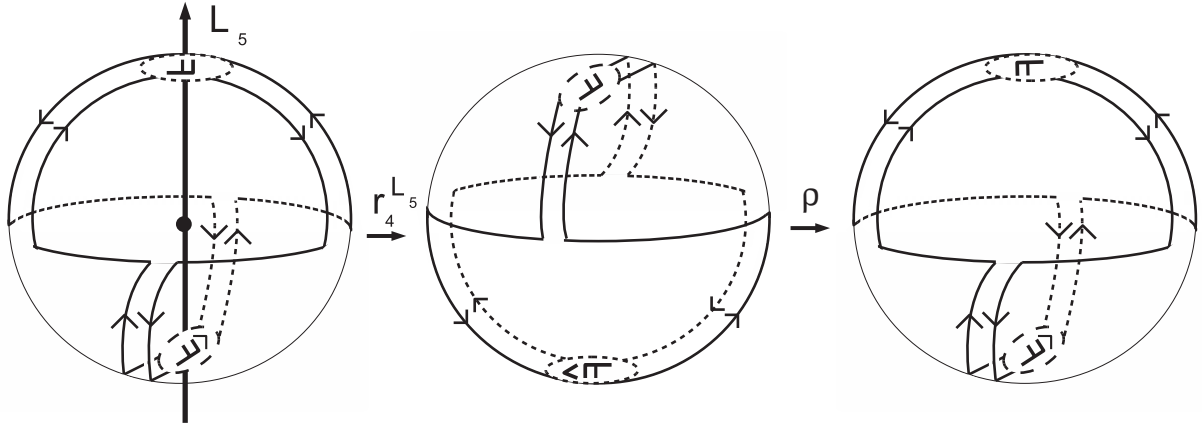}
\caption{The action of $R^{L_5}_4$}
\end{figure}

We see that  $\Pi$ is invariant by flypes and $R^{L_5}_4$ if and only if $F \sim F^*$.

{\bf Analysis of (v) $R^{L_5}_{-4}$}. The action of $R^{L_5}_{-4}$ is  very similar to that of $R^{L_5}_4$. Again, the conclusion is that $\Pi$ is invariant by flypes and $R^{L_5}_{-4}$ if and only if $F \sim F^*$.

{\bf Analysis of (vi) $R^{L_3}_2$ and of (vii)  $R^{L_4}_2$}. The action of $R^{L_3}_2$ is represented in Figure 19.

\begin{figure}[ht]    
   \centering
    \includegraphics[scale=0.35]{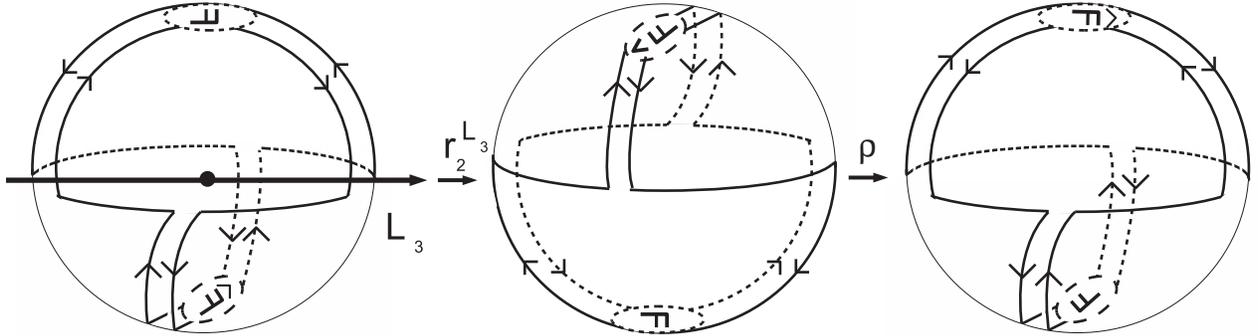}
\caption{The action of $R^{L_3}_2$}
\end{figure}

Here the situation is different from the preceding ones. We can transform $\Pi$ into $\widehat {\Pi}$ by flypes and $R^{L_3}_2$ or  $R^{L_4}_2$ if and only if $F \sim F^h \sim F^v$. Since $R^{L_3}_2$ and  $R^{L_4}_2$ reverse the orientation of $S^3$, these symmetries  are not pertinent to achirality.

{\bf Analysis of (viii) $r^{L_5}_2$}. Arguments quite similar to those above show that this symmetry is not suitable for achirality.

{\bf Summary.} 
\\
The symmetries $r^{L_1}_2$ and $r^{L_1}_2$ are relevant for -achirality. Note that these two symmetries are of order 2.
\\
The symmetries $R^{L_5}_4$ and $R^{L_5}_{-4}$ are relevant for +achirality. Note that these two symmetries are of order 4. 
\\
The other symmetries are not relevant to achirality.

{\bf Completing the proof of the proposition.} 
\\
The proposition states conditions on $F$ in order that $K$ be + or -achiral.
\\
That these conditions are sufficient is clear from the analysis of the figures. See also ``Remarks on the sufficiency of the conditions in Propositions 6.1 and 6.2" below.
\\
That these conditions are necessary is deep. Here, the fact that $\Pi$ is a minimal alternating projection is essential. The necessity relies on the following arguments:
\\
(a) the proof that an alternating $\pm$achiral knot must have a projection of Type I or II;
\\
(b) the ``important observation" that a projection of Type I or II is unique up to flypes in the two tangles;
\\
(c) the exhaustive analysis performed above of the symmetries of projections of Type I and II.

{\bf End of proof of the proposition.}

\subsection{Achirality of projections of Type II and completion of the proof
\\
 of Tait's Conjecture}

The careful analysis performed in the proof of Proposition 6.1 easily leads  to the next proposition.

\begin{proposition}
Let $K$ be an alternating knot which possesses a minimal projection  of Type II. Then:
\\
$K$ is always -achiral. Furthermore this symmetry is realised by a rotation of $S^2$ of angle $\pi$. 
\\
$K$ is +achiral if and only if $F \sim F^h \sim F^v$. In that case $K$ also admits a projection of Type I.
\end{proposition}

{\bf Remarks on the sufficiency of the conditions in Propositions 6.1 and 6.2}. The figures presented in the proof of Proposition 6.1 and the comments which accompany them show that the conditions expressed in the statements of  Propositions 6.1 and 6.2 are always sufficient. Explicitly:
\\
(1) Suppose that the knot $K$ possesses a (non-necessarily alternating) projection $\Pi$  of Type I. Then:
\\
$K$ is +achiral  if $\Pi$ satisfies $F \sim F^*$;
\\
$K$ is -achiral if $\Pi$ satisfies $F \sim F^h$ or $F \sim F^v$.
\\
(2) Suppose that the knot $K$ possesses a (non-necessarily alternating) projection $\Pi$  of Type II. Then:
\\
$K$ is +achiral if $\Pi$ satisfies $F \sim F^h \sim F^v$;
\\
$K$ is always -achiral.

The next theorem summarises the situation for projections with an invariant Haseman circle.

\begin{theorem}
Let $K$ be an alternating achiral knot and let $\Pi$ be  a minimal projection of $K$ with an invariant  Haseman circle. Then:
\\
(1) $K$ is +achiral if and only if $K$ possesses a minimal projection of Type I such that $F \sim F^*$;
\\
(2) $K$ is -achiral if and only if $K$ possesses a minimal projection of Type II. In this case, $\Pi$ can be transformed into its mirror image $\widehat {F}$ by a rotation of order 2. 
\end{theorem}

\begin{corollary}
Tait's Conjecture  for alternating -achiral knots is true.
\end{corollary}

{\bf Proof of the corollary.}
\\
At the end of Section 4, we  divided the proof of Tait's Conjecture into three cases named A, B and C. The proof of Case A is given in Section 5, which contains also a proof of Case B. The theorem above implies that Tait's Conjecture is true in Case B and C.
\\
{\bf End of proof of the corollary}

{\bf Comments.} 1) Roughly speaking, when a Haseman circle is invariant we can expect to have symmetries of order 2 (for -achirality) or  of order 4 (for +achirality). If we wish to have symmetries of higher order (a power of 2 and only in case of +achirality) we have to look for an invariant jewel. 
\\
2) Going back in time, we could say that Tait and Haseman discovered and studied achirality when a Haseman circle is invariant. Tait's efforts were devoted to -achirality (hence to symmetries of order 2) and Haseman's  were devoted to +achirality and  symmetries of order 4.

\section{The two checkerboard graphs}

In this section we show how the methods developed in this paper can be used to answer questions about checkerboards. There are two points we wish to emphasise. First, plain achirality is not precise enough. It is necessary to clearly distinguish between $+$ and $-$ achirality. It is also important to consider the nature of the fixed point of $\varphi : \mathcal {A} (K) \rightarrow  \mathcal {A} (K)$ as presented at the end of Section 4. We propose the following definition.

\begin{definition}
An achiral alternating knot $K$ is {\bf quasi-polyhedral} if the fixed point of $\varphi$ is either a jewel or a Haseman circle adjacent to two jewels (Case A or Case B).
\\
The knot $K$ is {\bf quasi-arborescent} if the fixed point of $\varphi$ is a Haseman circle adjacent to two band diagrams (Case C).
\end{definition}

Let $\Pi$ be a minimal projection of an alternating knot $K$. We denote by $G(\Pi)$ and $G^*(\Pi)$ the two checkerboard graphs associated to $\Pi$. For more information about these graphs see for instance \cite{QW1} where the graphs are denoted by $\Delta$ and $\Lambda$. 

The question whether $G(\Pi)$ and $G^*(\Pi)$ are isomorphic as planar graphs (i.e. as graphs embedded in the 2-sphere) is related to the achirality of $K$ as already noted and exploited by Tait. Louis Kauffman made the following conjecture. See \cite{KJ}.

{\bf Kauffman's Conjecture.} An alternating achiral knot has a minimal projection $\Pi$ such that $G(\Pi)$ is isomorphic to $G^*(\Pi)$.

To explicitly distinguish between $+$ and $-$ achirality it is necessary to make precise what ``isomorphic" means.

\begin{definition}
Two planar graphs $\Gamma$ and $\Gamma^*$ are {$\bf \pm$equivalent} if there exists a diffeomorphism $g : S^2 \rightarrow S^2$ of degree $\pm 1$ such that $g(\Gamma) = \Gamma^*$.
\end{definition} 

The next proposition is presumably well known.

\begin{proposition}
Let $\Pi$ be a minimal projection of an alternating knot $K$.
\\
i) One can move $\Pi$ to $-\widehat {\Pi}$ by an orientation  preserving diffeomorphism of $S^2$ if and only if $G(\Pi)$ is $+$equivalent to $G^*(\Pi)$.
\\
ii) One can move $\Pi$ to $+\widehat {\Pi}$ by an orientation  preserving diffeomorphism of $S^2$ if and only if $G(\Pi)$ is $-$equivalent to $G^*(\Pi)$.
\end{proposition}
{\bf Proof .}
For case i)($-$achirality) see  [15, Prop.6.2] and for case ii)($+$achirality) see [15, Prop.7.2].

The next theorem is a consequence of the proof of  Tait's  Conjecture.

\begin{theorem}
Kauffman's Conjecture is true for $-$achiral alternating knots. More precisely, every $-$achiral alternating knot has a minimal projection $\Pi$ such that $G(\Pi)$ is $+$equivalent to $G^*(\Pi)$.
\end{theorem}

However it is known that Kauffman's Conjecture is not true  in general . A counterexample was provided by Dasbach-Hougardy in \cite{DH}. From the above theorem, we deduce that counterexamples to Kauffman's Conjecture are necessarily $+$achiral and not $-$achiral.

In \cite{KJ} Kauffman and Jablan observe that the Dasbach-Hougardy counterexample is arborescent. It is easy to construct infinite families of such knots. Hence Kauffman-Jablan put forward the following conjecture.

{\bf Kauffman-Jablan's Conjecture.} Let $K$ be an alternating knot which is  a counterexample to Kauffman's Conjecture (such a knot is called Dasbach-Hougardy in \cite{KJ}). Then $K$ is arborescent.  

In a future publication \cite{EQ} we shall study in details $+$achiral alternating knots. Among other things, we shall prove that Kauffman-Jablan's Conjecture is not true. More precisely we have the following result.

\begin{theorem}
For every integer $a \geq 2$ there exists an alternating knot $K$ which is:
\\
1) $+$achiral of period $2^a$,
\\
2) quasi-polyhedral, 
\\
3) Dasbach-Hougardy.
\end{theorem}

Here is an example of such a knot where the tangle $F$ is pictured in Figure 14 while asking the projection to be alternating. See also the third comment at the end of Section 5.

\begin{figure}[ht]    
   \centering
    \includegraphics[scale=0.29]{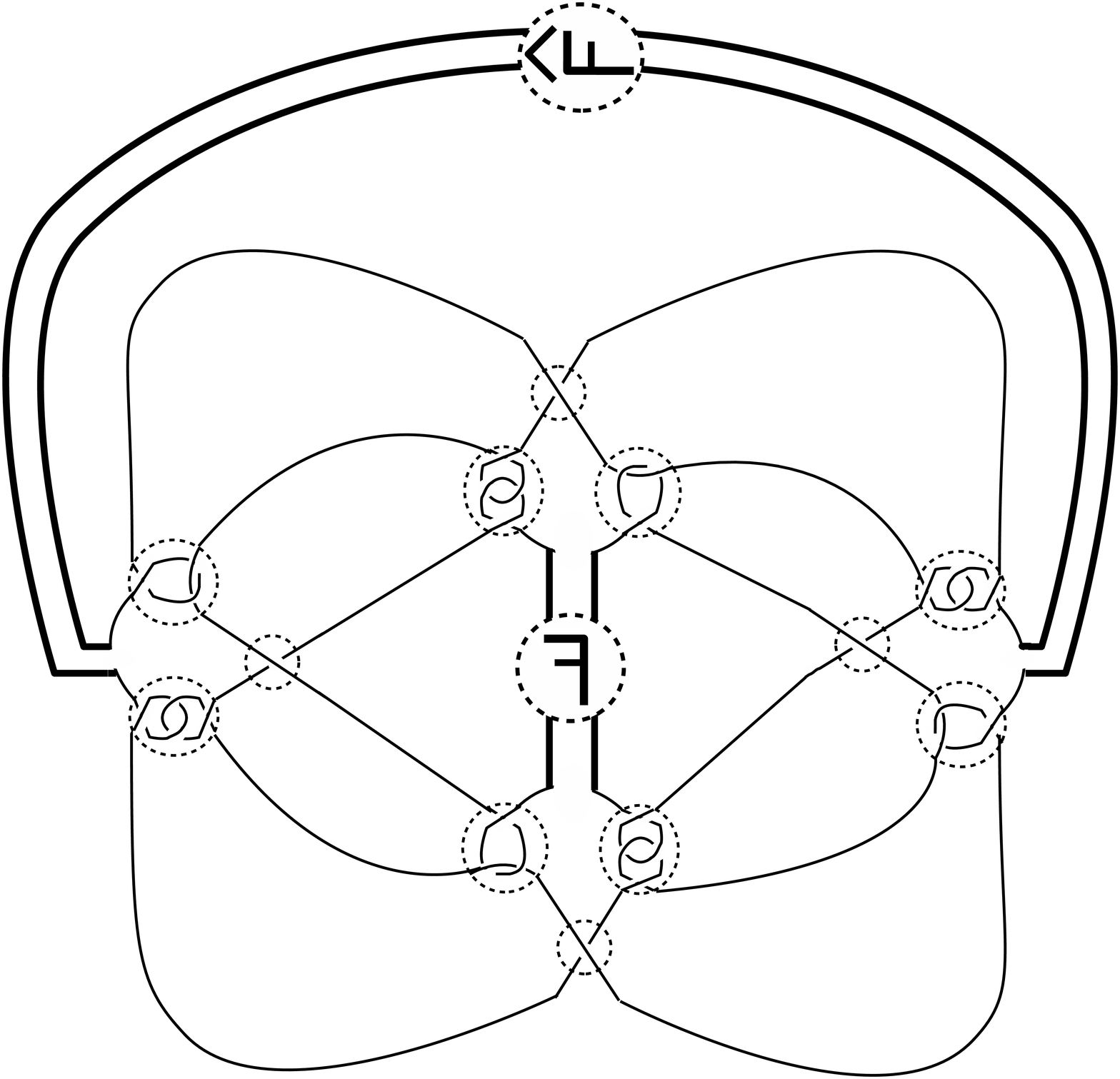}
\end{figure}

\section*{Acknowledgments} 
We  thank Le Fonds National Suisse de la Recherche Scientifique for its  support.
\\
We wish to thank for John Steinig for his generous help.

\newpage


 {\bf Address.}

Section de math\'ematiques
\\
Universit\'e de Gen\`eve
\\
Cp 64
\\
CH-1211 GENEVE 4
\\
SWITERLAND

nicola.ermotti@unige.ch
\\
  cam.quach@math.unige.ch
\\
claude.weber@math.unige.ch

\end{document}